\documentclass{amsproc}
\usepackage{amsmath, amscd, amsthm, amsfonts, amssymb, mathrsfs}

\newtheorem{theorem}{Theorem}[section]

\theoremstyle{definition}

\newtheorem{problem}{Problem}

\newtheorem{conjecture}{Conjecture}
\newtheorem{corollary}{Corollary}

\theoremstyle{remark}

\input xy
\xyoption{all}

\numberwithin{equation}{section}

\DeclareMathOperator{\Ree}{R}
\DeclareMathOperator{\Reeb}{^2G_2}
\DeclareMathOperator{\Orth}{O}

\newcommand{\floor}[1]{\lfloor#1\rfloor}

\newcommand{\Fk}{\Bbb{F}_{2^k}}

\newcommand{\Sz}{\mathrm{Sz}}

\newcommand{\Sp}{\mathrm{Sp}}

\newcommand{\PSL}{\mathrm{PSL}} 
\newcommand{\PGL}{\mathrm{PGL}}

\newcommand{\PGaL}{\mathrm{P\Gamma L}} 
\newcommand{\PSiL}{\mathrm{P\Sigma L}}

\newcommand{\Sym}{\mathrm{S}}
\newcommand{\Alt}{\mathrm{A}}

\newcommand{\Cyc}{\mathrm{C}}

\newcommand{\rz}{\rho_0}

\newcommand{\rd}{\rho_2}
\newcommand{\lb}{\lbrace}
\newcommand{\rb}{\rbrace}

\begin{document}

\title{String C-group representations of almost simple groups: a survey}

\author{Dimitri Leemans}
\address{Dimitri Leemans, Universit\'{e} Libre de Bruxelles,
D\'{e}partement de Math\'{e}matique,
C.P.216 Alg\`ebre et Combinatoire,
Bld du Triomphe, 1050 Bruxelles,
Belgium}
\email{dleemans@ulb.ac.be}

%    General info
\subjclass{Primary 20D06, 52B11}
\date{}


\keywords{String C-group representations, abstract regular polytopes, almost simple groups}

\begin{abstract}
This survey paper aims at giving the state of the art in the study of string C-group representations of almost simple groups. It also suggest a series of problems and conjectures to the interested reader.
\end{abstract}

\maketitle

\tableofcontents

\section{Introduction}\label{intro}
Polytopes, and in particular  regular polytopes, have been studied by mathematicians for millenia. The recent monograph of McMullen and Schulte~\cite{ARP} came as the first comprehensive up-to-date book on abstract regular polytopes after more than twenty years of rapid development and it describes a rich new theory that benefits from an interplay between several areas of mathematics including geometry, algebra, combinatorics, group theory and topology. 
In recent years, abstract regular polytopes whose automorphism groups are (almost) simple groups have been studied extensively. This survey paper aims to give the state of the art on the subject, and to give paths for further research in that field.

The project surveyed in this paper started in 2003 when Michael Hartley contacted us about an abstract regular polytope he had found, of type $\{5,3,5\}$, whose automorphism group is the first Janko group $J_1$, one of the 26 sporadic simple groups. It turned out we knew the existence of that polytope since 1999 when we published an atlas of regular thin geometries for small groups~\cite{Lee97c}.
Unfortunately (or in fact fortunately), that polytope was not mentioned in this atlas as, at the time, we had computed all thin regular residually connected incidence geometries whose automorphism group is $J_1$ but there were too many of them (almost 3000) to list them in the paper.
That polytope was the start of a very fruitful collaboration with Hartley and changed our research life as we then decided to focus on polytopes.

Hartley was in fact willing to determine whether the universal locally projective polytope of type $\{5,3,5\}$ was finite or infinite. That was the missing piece to finalise the classification of universal locally projective polytopes.
We eventually found out that it is finite and its automorphism group is isomorphic to $J_1\times PSL(2,19)$, a very surprising result~\cite{HL2003}.
Indeed, Gr\"unbaum found in 1977 a rank four polytope which he called the 11-cell~\cite{grunbaum}. The automorphism group of the 11-cell is the group $PSL(2,11)$ and Gr\"unbaum obtained it by taking the Coxeter group of type $[3,5,3]$ and adding two relations, to quotient the facets and vertex-figures and make them hemi-icosahedra and hemi-dodecahedra respectively. In the same spirit, Coxeter found the 57-cell~\cite{cox57} by adding two relations to the Coxeter group of type $[5,3,5]$ to make the 
facets and vertex-figures hemi-dodecahedra and hemi-icosahedra respectively.
It turns out that in the case of the 11-cell, adding only one of the two relations implies the other relation, but in the case of the 57-cell, it is not the case: adding only one relation gives the group $J_1\times PSL(2,19)$ as we found out with Hartley.

We then met in Brussels for a couple of weeks of intense work and decided at the time it would be great to build atlases of polytopes. Hartley focused on polytopes with a fixed number of flags~\cite{Halg} while we focused on polytopes whose automorphism groups are almost simple groups with Laurence Vauthier~\cite{LV}. We started collecting a sensible amount of computer-generated data on the subject and used them to state conjectures.

The aim of this paper is to give a survey of what has been done on that subject over the last fifteen years, to also give some results that were obtained on C-groups, to state some conjectures and give open problems. 

The three big questions we have tried to answer for families of groups over the years are the following ones. Given a group $G$ (eventually belonging to a family of groups as, for instance finite simple groups, or alternating groups),

\begin{enumerate}
\item determine the possible ranks of string C-group representations for $G$;
\item in particular, determine the highest rank of a string C-group representation of $G$;
\item enumerate all string C-group representations for $G$;
\end{enumerate}
These questions make sense in trying to find nice geometric structures on which these groups act as automorphism groups, but they also tell us something about quotients of Coxeter groups -- a string C-group representation being a smooth quotient of a Coxeter group -- as well as independent generating sets of the groups $G$.

One can consider these same three questions for chiral polytopes. In particular, obtaining enumeration data for both types of polytopes would shed light on their relative abundance. However, we constrain this survey to questions about regular polytopes.
Observe that the enumeration question is usually very hard to answer for families of groups and that, apart for shedding light on the relative abundance of regular polytopes versus chiral polytopes, its interest might be considered lower than the interest of the first two questions.

Even though almost simple groups appear scarcely among the set of all groups, trying to answer questions (1), (2) and (3) for the finite simple groups (and their automorphism groups) seems natural in order to get a better geometric understanding of some of them. Moreover, for finite simple groups, involutions have been thoroughly studied in the process of their classification.
The families of (almost) simple groups not mentioned in this survey are of course also very interesting to look at but, currently, no results are known for them. So any discovery for them is most welcome.

\section{Preliminaries}\label{preliminaries}
Let us start by defining the concepts and fixing the notation needed to understand this survey.
 
\subsection{C-groups and string C-groups}
As it is well known that abstract regular polytopes and string C-groups are in one-to-one correspondence (see for instance~\cite[Section 2E]{ARP}), and since it is much easier to define string C-groups, we frame our discussion in the language of string C-groups.

Let $G$ be a group and $X$ a set of involutions of $G$.
If $\langle X \rangle = G$, we say that $(G;X)$ is a {\em group generated by involutions} or {\em ggi} for short.
If there is an ordering of $X$, say 
 $\{\rho_0,\ldots, \rho_{n-1}\}$ where $\rho_i \rho_j = \rho_j\rho_i$ for every $i,j\in \{0,\ldots, n-1\}$ such that $|i-j|>1$, then $(G;\{\rho_0,\ldots, \rho_{n-1}\})$ is 
a {\em string group generated by involutions} or {\em sggi} for short.

When $(G;\{\rho_0,\ldots, \rho_{n-1}\})$ is an sggi, we assume, without loss of generality, that the involutions are ordered in such a way that $\forall i,j\in\{0,\ldots, r-1\}$, if $\;|i-j|>1$ then $(\rho_i\rho_j)^2=1$ (this property is called the {\em commuting property}).

Let $G$ be a group and $\{\rho_0,\ldots, \rho_{n-1}\}$ be a set of involutions of $G$.
The pair $(G;\{\rho_0,\ldots, \rho_{n-1}\})$ satisfies the {\em intersection property} if
\[\forall I,J \subseteq \{0,\ldots,n-1\}, \; \langle \rho_i \mid i \in I\rangle \cap \langle \rho_j \mid j 
 \in J\rangle = \langle \rho_k \mid k \in I\cap J\rangle.\]
 
If a ggi $(G;\{\rho_0,\ldots, \rho_{n-1}\})$ satisfies the intersection property, it is called a {\em C-group representation of $G$} (or C-group for short).
A C-group that satisfies the commuting property is called a {\em string C-group}. 

The integer $n$ is the {\em rank} of the (string) C-group $(G;\{\rho_0,\ldots, \rho_{n-1}\})$.

To a string C-group $(G;\{\rho_0,\ldots, \rho_{n-1}\})$, we can associate a {\em Schl\"afli type}, that is a sequence $\{p_1, \ldots, p_{n-1}\}$ where $p_i = o(\rho_{i-1}\rho_i)$. Observe that if one of the $p_i$'s is equal to 2, then the group $G$ is a direct product of two non-trivial subgroups. If that is the case, we say that $(G;\{\rho_0,\ldots, \rho_{n-1}\})$ is {\em reducible}; Otherwise $(G;\{\rho_0,\ldots, \rho_{n-1}\})$ is called {\em irreducible}. Obviously, if $G$ is simple and $(G;\{\rho_0,\ldots, \rho_{n-1}\})$ is a string C-group, then $(G;\{\rho_0,\ldots, \rho_{n-1}\})$ is irreducible.
We take as convention that a symbol $p_i$ appearing $k$ times in adjacent places in a Schl\"afli type can be replaced by $p_i^k$ instead of writing $k$ times $p_i$. So for instance the Schl\"afli type of the 4-simplex can be written as $\{3,3,3\}$ or as $\{3^3\}$.  

Given a group $G$ and two sets $S$ and $T$ of involutions of $G$ such that $(G,S)$ and $(G,T)$ are  C-groups, we say that $(G,S)$ and $(G,T)$ are {\em isomorphic} if there exists an element $g\in Aut(G)$ such that $g$ maps $S$ onto $T$.

The {\em (string) C-rank} of a group is the largest possible rank of a (string) C-group representation for that group.

\subsection{Permutation representation graphs and CPR graphs}
Let $G$ be a group of permutations acting on a set $\{1,\,\ldots,\,n\}$.
Let $S:=\{\rho_0,\ldots,\rho_{r-1}\}$ be a set of $r$ involutions of $G$ such that $G=\langle S \rangle$.
We define the \emph{permutation representation graph} $\mathcal{G}$  of $G$ as the $r$-edge-labeled multigraph with vertex set $V(\mathcal{G}) := \{1,\,\ldots,\,n\}$. The edge-set of $\mathcal{G}$ is the set $\{\{a,a\rho_i\} : a \in V(\mathcal{G}), i\in \{0, \ldots, r-1\}| a\neq a\rho_i\}$ where every edge $\{a,a\rho_i\}$ has label $i$.
When $(G,S)$ is a string C-group representation, the permutation representation graph $\mathcal{G}$ is called a \emph{CPR graph}, as defined in \cite{pcpr}.

\section{Simple groups and rank three string C-group representations}\label{nuzhin}
In 1980, it was asked in the Kourovka Notebook (Problem 7.30) which finite simple groups can be generated by three involutions, two of which commute.
This problem was solved by Nuzhin and Mazurov in \cite{Mazurov03,n1,n2,n3,n4}. The groups $PSU_4(3)$ and $PSU_5(2)$, although mentioned by Nuzhin as being generated by three involutions, two of which commute, have recently been discovered not to have such generating sets by Martin Macaj and Gareth Jones (personal communication of Jones, checked independently using {\sc Magma}). We summarize below the accurate solution to Problem 7.30.

\begin{theorem}[Nuzhin - Mazurov - Macaj - Jones]\label{nmmj}
Every non-abelian finite simple group can be generated by three involutions, two of which commute, with the following exceptions:
\[
\begin{array}{l}
PSL_3(q),\,PSU_3(q),\,PSL_4(2^n), \,PSU_4(2^{n}),
A_6,\\ 
A_7, PSU_4(3), \,PSU_5(2), 
\,M_{11},\, M_{22}, \,M_{23}, \,McL.
\end{array}
\]
\end{theorem}

If $G$ is a simple group generated by three involutions, two of which commute, $G$ together with these three involutions form a string C-group representation by the following result due to Marston Conder and Deborah Oliveros.
\begin{theorem}~\cite[Corollary 4.2]{CO2013}\label{CO}
If $G$ is a finite non-abelian simple group, or more generally any finite group with no non-trivial cyclic normal subgroup, then every smooth homomorphism from the $[k,m]$ Coxeter group onto $G$ gives rise to a string C-group representation of rank three for $G$.
\end{theorem}
It turns out that Theorem~\ref{nmmj} holds when removing the hypothesis on the rank as proven by Adrien Vandenschrick in~\cite{Adrien}: the exceptions in rank three do not have string C-group representations of higher rank.

\section{Symmetric and alternating groups}

\subsection{String C-group representations of $S_n$}
Symmetric groups gained our attention early in this research project. They were among the only ones in the data we collected that gave string C-group representations of large rank.
It was known for a long time~\cite{Moo1896} that the group $S_n$ has a string C-group representation of rank $n-1$, namely $(S_n;\{(1,2),(2,3),\ldots, (n-1,n)\})$.
A recent result of Julius Whiston~\cite{Whis00}, showing that the largest size of a set of independent generators of $S_n$ is $n-1$, implies that the string C-rank of $S_n$ is $n-1$.

Sjerve and Cherkassoff showed in~\cite{SC94} that $S_n$ is a group generated by three involutions, two of which commute, provided that $n\geq 4$.
Their examples satisfy the intersection property and therefore are rank three string C-group representations by Theorem~\ref{CO}. 
%Observe that there are also triples of involutions, none of which commuting, that generate $S_n$ and satisfy the intersection property.

\begin{theorem}\cite[Theorem 1.2]{SC94}
Every group $S_n$ with $n\geq 4$ has at least one string C-group representation of rank three.
\end{theorem}
Earlier work by Conder~\cite{Con80,Con81} covers all but a few cases of the results of~\cite{SC94} for the symmetric groups.
As Conder pointed out to us, these days, it takes a few seconds to handle the missing cases for $S_n$ with {\sc Magma}~\cite{BCP97}.

Together with Maria Elisa Fernandes, we determined the number of non-isomorphic string C-group representations of $S_n$ of rank $n-1$ and $n-2$~\cite{fl, Corr}. Then with Fernandes and Mark Mixer we also determined the representations of rank $n-3$ and $n-4$~\cite{extension}.

\begin{theorem}\label{maintheorem}\cite[Theorem 1.1]{extension}
Let $1\leq i \leq 4$, and $n\geq 3+2i$ when $r = n-i$.
If $\Gamma$ is a string C-group representation of a group $G$, of rank $r\geq n-i$ with a connected diagram and $G$ is isomorphic to a transitive group of degree $n$  then $G$ or its dual is isomorphic to $S_n$ and the CPR graph is one of those listed in Table~\ref{7cases}.
\end{theorem}

\begin{small}
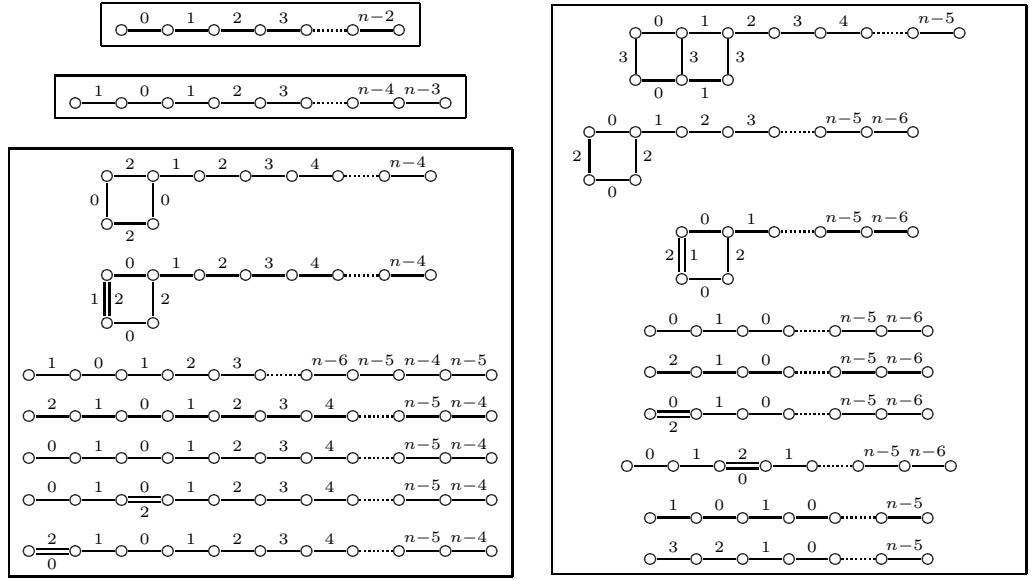
\begin{figure}[h]
\begin{tabular}{cc}
\begin{tabular}{c}
\begin{tabular}{|c|}
\hline
$ \xymatrix@-1pc{ *+[o][F]{}    \ar@{-}[r]^0    &  *+[o][F]{}   \ar@{-}[r]^1  & *+[o][F]{}  \ar@{-}[r]^2 & *+[o][F]{}  \ar@{-}[r]^3 & *+[o][F]{} \ar@{.}[r] &  *+[o][F]{}  \ar@{-}[r]^{n-2} & *+[o][F]{}}$\\
\hline
\end{tabular}\\

\\

\begin{tabular}{|c|}
\hline
$ \xymatrix@-1pc{ *+[o][F]{}   \ar@{-}[r]^1  & *+[o][F]{}  \ar@{-}[r]^0    &  *+[o][F]{}   \ar@{-}[r]^1  & *+[o][F]{}  \ar@{-}[r]^2 & *+[o][F]{}  \ar@{-}[r]^3 & *+[o][F]{} \ar@{.}[r] &   *+[o][F]{}  \ar@{-}[r]^{n-4} & *+[o][F]{}  \ar@{-}[r]^{n-3} & *+[o][F]{} }$\\
\hline
\end{tabular}\\
\\
\begin{tabular}{|c|}
\hline
$ \xymatrix@-1pc{
     *+[o][F]{}   \ar@{-}[r]^2  \ar@{-}[d]_0     & *+[o][F]{}  \ar@{-}[r]^1 \ar@{-}[d]^0  & *+[o][F]{}  \ar@{-}[r]^2 & *+[o][F]{}  \ar@{-}[r]^3 & *+[o][F]{}  \ar@{-}[r]^4 & *+[o][F]{}  \ar@{..}[r] & *+[o][F]{}  \ar@{-}[r]^{n-4}& *+[o][F]{}  \\
     *+[o][F]{}  \ar@{-}[r]_2     & *+[o][F]{}    &  & & & & &  }$
\\
$ \xymatrix@-1pc{
     *+[o][F]{}   \ar@{-}[r]^0  \ar@{=}[d]_1^2     & *+[o][F]{}  \ar@{-}[r]^1 \ar@{-}[d]^2  & *+[o][F]{} \ar@{-}[r]^2 & *+[o][F]{}  \ar@{-}[r]^3 & *+[o][F]{}  \ar@{-}[r]^4 & *+[o][F]{}  \ar@{..}[r] & *+[o][F]{}  \ar@{-}[r]^{n-4}& *+[o][F]{}  \\
     *+[o][F]{}  \ar@{-}[r]_0     & *+[o][F]{}    &  & & & & &  }$
     \\
$ \xymatrix@-1pc{ *+[o][F]{}   \ar@{-}[r]^1  & *+[o][F]{}  \ar@{-}[r]^0    &  *+[o][F]{}   \ar@{-}[r]^1  & *+[o][F]{}  \ar@{-}[r]^2 & *+[o][F]{}  \ar@{-}[r]^3 & *+[o][F]{} \ar@{.}[r] &  *+[o][F]{}  \ar@{-}[r]^{n-6} & *+[o][F]{}  \ar@{-}[r]^{n-5} & *+[o][F]{}  \ar@{-}[r]^{n-4} & *+[o][F]{}  \ar@{-}[r]^{n-5} & *+[o][F]{}}$
\\
 $ \xymatrix@-1pc{ *+[o][F]{}   \ar@{-}[r]^2  & *+[o][F]{}  \ar@{-}[r]^1    &  *+[o][F]{}   \ar@{-}[r]^0  & *+[o][F]{}  \ar@{-}[r]^1 & *+[o][F]{}  \ar@{-}[r]^2 & *+[o][F]{}  \ar@{-}[r]^3  & *+[o][F]{}  \ar@{-}[r]^4  &*+[o][F]{} \ar@{.}[r] &  *+[o][F]{}  \ar@{-}[r]^{n-5} & *+[o][F]{}  \ar@{-}[r]^{n-4} & *+[o][F]{}  }$
 \\
$ \xymatrix@-1pc{ *+[o][F]{}   \ar@{-}[r]^0  & *+[o][F]{} \ar@{-}[r]^1  &  *+[o][F]{} \ar@{-}[r]^0 & *+[o][F]{} \ar@{-}[r]^1& *+[o][F]{} \ar@{-}[r]^2  & *+[o][F]{}  \ar@{-}[r]^3  & *+[o][F]{}  \ar@{-}[r]^4 & *+[o][F]{} \ar@{.}[r] & *+[o][F]{}  \ar@{-}[r]^{n-5} & *+[o][F]{}  \ar@{-}[r]^{n-4} & *+[o][F]{}}$
\\
 $ \xymatrix@-1pc{ *+[o][F]{}   \ar@{-}[r]^0  & *+[o][F]{} \ar@{-}[r]^1  &  *+[o][F]{} \ar@{=}[r]^0_2 & *+[o][F]{} \ar@{-}[r]^1& *+[o][F]{} \ar@{-}[r]^2  & *+[o][F]{}  \ar@{-}[r]^3  & *+[o][F]{}  \ar@{-}[r]^4 & *+[o][F]{} \ar@{.}[r] & *+[o][F]{}  \ar@{-}[r]^{n-5} & *+[o][F]{}  \ar@{-}[r]^{n-4} & *+[o][F]{}}$
 \\
$ \xymatrix@-1pc{ *+[o][F]{}   \ar@{=}[r]^2_0  & *+[o][F]{} \ar@{-}[r]^1  &  *+[o][F]{} \ar@{-}[r]^0 & *+[o][F]{} \ar@{-}[r]^1& *+[o][F]{} \ar@{-}[r]^2  & *+[o][F]{}  \ar@{-}[r]^3  & *+[o][F]{}  \ar@{-}[r]^4 & *+[o][F]{} \ar@{.}[r] & *+[o][F]{}  \ar@{-}[r]^{n-5} & *+[o][F]{}  \ar@{-}[r]^{n-4} & *+[o][F]{}}$  \\
\hline
\end{tabular}
\end{tabular}
&
\begin{tabular}{|c|}
\hline
$\xymatrix@-1pc{
         *+[o][F]{}   \ar@{-}[d]_3  \ar@{-}[r]^0     & *+[o][F]{}  \ar@{-}[d]^3  \ar@{-}[r]^1  & *+[o][F]{} \ar@{-}[d]^3 \ar@{-}[r]^2 & *+[o][F]{} \ar@{-}[r]^3 & *+[o][F]{} \ar@{-}[r]^4& *+[o][F]{} \ar@{.}[r] &*+[o][F]{} \ar@{-}[r]^{n-5}&*+[o][F]{}\\
      *+[o][F]{}  \ar@{-}[r]_0     & *+[o][F]{}  \ar@{-}[r]_1   &  *+[o][F]{}&  &&&&}$
\\
$\xymatrix@-1pc{
      *+[o][F]{}   \ar@{-}[d]_2  \ar@{-}[r]^0     & *+[o][F]{}  \ar@{-}[d]^2  \ar@{-}[r]^1  & *+[o][F]{}\ar@{-}[r]^2 & *+[o][F]{}\ar@{-}[r]^3 & *+[o][F]{}\ar@{.}[r]& *+[o][F]{}\ar@{-}[r]^{n-5}&*+[o][F]{}\ar@{-}[r]^{n-6} &*+[o][F]{}\\
      *+[o][F]{}  \ar@{-}[r]_0     & *+[o][F]{}    &   &&&&&&&}$
\\
$\xymatrix@-1pc{
     *+[o][F]{}   \ar@{-}[r]^0  \ar@{=}[d]_2^1     & *+[o][F]{}  \ar@{-}[r]^1 \ar@{-}[d]^2 & *+[o][F]{}  \ar@{.}[r] &*+[o][F]{}  \ar@{-}[r]^{n-5} &*+[o][F]{}  \ar@{-}[r]^{n-6} &*+[o][F]{}\\
     *+[o][F]{}  \ar@{-}[r]_0     & *+[o][F]{}    & &&&   }$
\\
$\xymatrix@-1pc{ *+[o][F]{}   \ar@{-}[r]^0 & *+[o][F]{} \ar@{-}[r]^1  &  *+[o][F]{} \ar@{-}[r]^0 & *+[o][F]{} \ar@{.}[r] & *+[o][F]{}  \ar@{-}[r]^{n-5}& *+[o][F]{}  \ar@{-}[r]^{n-6}  & *+[o][F]{}}$
\\
$\xymatrix@-1pc{ *+[o][F]{}   \ar@{-}[r]^2  & *+[o][F]{} \ar@{-}[r]^1  &  *+[o][F]{} \ar@{-}[r]^0 & *+[o][F]{} \ar@{.}[r] & *+[o][F]{}  \ar@{-}[r]^{n-5}& *+[o][F]{}  \ar@{-}[r]^{n-6}  & *+[o][F]{}}$
\\
$\xymatrix@-1pc{ *+[o][F]{}   \ar@{=}[r]^0_2  & *+[o][F]{} \ar@{-}[r]^1  &  *+[o][F]{} \ar@{-}[r]^0 & *+[o][F]{} \ar@{.}[r] & *+[o][F]{}  \ar@{-}[r]^{n-5}& *+[o][F]{}  \ar@{-}[r]^{n-6}  & *+[o][F]{}}$
\\
$\xymatrix@-1pc{ *+[o][F]{}  \ar@{-}[r]^0& *+[o][F]{}  \ar@{-}[r]^{1} & *+[o][F]{}  \ar@{=}[r]^{2}_0 & *+[o][F]{}  \ar@{-}[r]^{1}  & *+[o][F]{} \ar@{.}[r]& *+[o][F]{}  \ar@{-}[r]^{n-5}  & *+[o][F]{}  \ar@{-}[r]^{n-6} & *+[o][F]{}}$
\\
$\xymatrix@-1pc{ *+[o][F]{}   \ar@{-}[r]^1  &*+[o][F]{}   \ar@{-}[r]^0  & *+[o][F]{} \ar@{-}[r]^1  &  *+[o][F]{} \ar@{-}[r]^0  & *+[o][F]{} \ar@{.}[r] &  *+[o][F]{}  \ar@{-}[r]^{n-5} & *+[o][F]{}}$
\\
$\xymatrix@-1pc{ *+[o][F]{}   \ar@{-}[r]^3 &*+[o][F]{}   \ar@{-}[r]^2  & *+[o][F]{} \ar@{-}[r]^1  &  *+[o][F]{} \ar@{-}[r]^0& *+[o][F]{}\ar@{.}[r] & *+[o][F]{}  \ar@{-}[r]^{n-5} & *+[o][F]{}}$
\\
\hline
\end{tabular}
\end{tabular}
\caption{CPR graphs of string C-group representations of rank $r\geq n-4$ for $S_n$.}\label{7cases}
\end{figure}
\end{small}

Also with Fernandes, we showed that there were no gaps in the set of possible ranks of string C-group representations of $S_n$.

\begin{theorem}\label{1}\cite[Theorem 3]{fl}
Let $n\geq 4$.
For every $r \in \{3, \ldots, n-1\}$, there exists at least one string C-group  representation of rank $r$ of $S_n$ for every $3\leq r \leq n-1$.
Its Schl\"afli type is $\{n-r+2, 6, 3^{r-3}\}$.
\end{theorem}

Table~\ref{sn} gives, for $S_n$ ($5\leq n\leq 14$), the number of pairwise nonisomorphic string C-group representations of rank $r$ ($3\leq r\leq n-1$). It suggests the following problem.
\begin{table}
\begin{tabular}{|c|c|c|c|c|c|c|c|c|c|c|c|c||}
\hline
$G$&Rk 3&Rk 4&Rk 5&Rk 6&Rk 7&Rk 8&Rk 9&Rk 10&Rk 11&Rk 12&Rk 13\\
\hline 
$S_5$&4&1&0&0&0&0&0&0&0&0&0\\
$S_6$&2&4&1&0&0&0&0&0&0&0&0\\
$S_7$&35&7&1&1&0&0&0&0&0&0&0\\
$S_8$&68&36&11&1&1&0&0&0&0&0&0\\
$S_9$&129&37&7&7&1&1&0&0&0&0&0\\
$S_{10}$&413 &203 &52 &13 &7&1&1&0&0&0&0\\
$S_{11}$&1221 &189 &43 &25 &9&7&1&1&0&0&0\\
$S_{12}$&3346 &940 &183 &75 &40 &9&7&1&1&0&0\\
$S_{13}$&7163 &863 &171 &123 &41 &35 &9&7&1&1&0\\
$S_{14}$&23126 &3945 &978 &303 &163 &54 &35 &9&7&1&1\\
\hline
\end{tabular}
\caption{The number of pariwise nonisomorphic string C-group representations of $S_n$ ($5\leq n\leq 14$).}\label{sn}
\end{table}
\begin{problem}\label{snenum}
Enumerate the string C-group representations of rank $r$ of $S_n$ for every $3\leq r \leq n-1$.
\end{problem}
This problem is solved for the four highest values of $r$ and $n$ large enough as pointed out in Theorem~\ref{maintheorem}. It suggests the following conjecture.

\begin{conjecture}\label{conjsn}
The number of string C-group representations of rank $n-i$ for $S_n$ with $1\leq i\leq (n-3)/2$ is a constant independent of $n$.
\end{conjecture}
This conjecture suggests the existence of a new integer sequence (not appearing in the Encyclopedia of Integer Sequences) starting with 1, 1, 7, 9 and whose next number is likely 35 as the experimental data of Table~\ref{sn} suggest.

A first step in trying to solve Problem~\ref{snenum} has been made by Kiefer and the author. We 
computed the number of unordered pairs of commuting involutions in $S_n$ up to conjugacy in $Aut(S_n)$ in~\cite{KL2013}. We obtained the following result.

\begin{theorem}\cite[Theorem 1.1]{KL2013}\label{pairsym}
Let $n>1$ be a positive integer.
Define $\lambda(k)$ and $\psi(k,n)$ as follows.
\begin{equation*}\begin{aligned}
\lambda(k) =\Big\lfloor{\Big(\frac{k}{2}+1\Big)^2}\Big\rfloor
\end{aligned}\end{equation*}
\begin{equation*}\begin{aligned}
\psi(k,n) & = \begin{cases}
\left[ \frac{1}{2} \left( 2k - n \right)\right]^2 + \frac{1}{2} \left( 2k - n \right) & \textrm{ if } n \textrm{ is even}, \\
 & \\
\left[ \frac{1}{2} \left( 2k - n -1\right)\right]^2 + 2k - n  & \textrm{ if } n \textrm{ is odd.} \\
\end{cases}
\end{aligned}\end{equation*}

There are, up to isomorphism,
\begin{equation*}\begin{aligned}
-\frac{3n}{2}+\sum_{k=1}^{n} \lambda(k)\cdot \Big(\frac{n-k}{2}+ 1\Big) - \frac{1}{2} \cdot \sum_{k= \floor{\frac{n}{2}}+1}^{n} \psi(k,n)
\end{aligned}\end{equation*}

unordered pairs of commuting involutions in $S_{2n}$ and $S_{2n+1}$ except for $S_6$ in which there are, up to isomorphism, five 
unordered pairs of commuting involutions.
\end{theorem}
A similar result was obtained for the alternating groups (see Theorem~\ref{pairalt}).
Table~\ref{snancommut} gives the number of pairwise nonisomorphic pairs of commuting involutions for $S_n$ (and $A_n$) for some values of $n$.
\begin{table}
\begin{tabular}{|c|c|c|}
\hline
{$n$}&{$\lb \rz, \rd \rb$, with $\rz,\rd \in \Sym(n)$}&{$\lb \rz, \rd \rb$, with $\rz,\rd \in \Alt(n)$}\\
\hline
\hline
1,2,3	&0	&0\\
4,5	&3	&1\\
6	&5	&1\\
7	&9	&2\\
8,9	&21	&7\\
10,11	&39	&10\\
12,13	&67	&21\\
14,15	&105	&28\\
16,17	&158	&48\\
18,19	&226	&61\\
20	&315	&93\\
30	&1169	&315\\
40	&3105	&855\\
50	&6774	&1795\\
\hline
\end{tabular}
\caption{Number of unordered pairs of commuting involutions in $S_n$ and $A_n$, up to conjugacy in $Aut(S_n)$.}\label{snancommut}
\end{table}
As the reader can see, the numbers appearing in this table are not very encouraging to solve Problem~\ref{snenum} for the cases not covered by Conjecture~\ref{conjsn}.

\subsection{String C-group representations of $A_n$}
Alternating groups were investigated in the same vein as symmetric groups.

Sjerve and Cherkassoff showed in~\cite{SC94} that $A_n$ is a group generated by three involutions, two of which commute, provided that $n\geq 4$ and $n\neq 6, 7$ or $8$.
Their examples satisfy the intersection property and therefore are rank three string C-group representations.

\begin{theorem}\cite[Theorem 1.1]{SC94}
Every group $A_n$ with $n=5$ or $n\geq 9$ has at least one string C-group representation of rank three.
\end{theorem}
Again, earlier work by Conder~\cite{Con80,Con81} covers all but a few cases of the results of~\cite{SC94} for the alternating groups.

When we started working on string C-group representations for the alternating groups with Fernandes and Mixer, we first collected experimental data for $A_n$ with $n\leq 12$.
Table~\ref{an} gives the number of pairwise nonisomorphic string C-group representations for $A_n$ ($5\leq n \leq 15$).
\begin{table}
\begin{center}
\begin{tabular}{||c|c|c|c|c|c|c||}
\hline
$G$&Rank 3&Rank 4&Rank 5&Rank 6&Rank 7&Rank 8\\
\hline
$A_5$&2&0&0&0&0&0\\
$A_6$&0&0&0&0&0&0\\
$A_7$&0&0&0&0&0&0\\
$A_8$&0&0&0&0&0&0\\
$A_9$&41&6&0&0&0&0\\
$A_{10}$&94&2&4&0&0&0\\
$A_{11}$&64&0&0&3&0&0\\
$A_{12}$&194&90&22&0&0&0\\
$A_{13}$&1558&102&25&10&0&0\\
$A_{14}$&4347&128&45&9&0&0\\
$A_{15}$&5820&158&20&42&6&0\\
\hline
\end{tabular}
\caption{The number of pariwise nonisomorphic string C-group representations of $A_n$ ($5\leq n\leq 15$).}\label{an}
\end{center}
\end{table}
A striking observation came for $A_{11}$. It was the first time we found a group that had string C-group representations whose set of ranks is not an interval (in this case, the set of ranks is $\{3,6\}$).

We proved that for each rank $r\geq 4$, there is at least one group $A_n$ that has a string C-group representation of rank $r$.

\begin{theorem}\cite[Theorem 1.1]{flm}
For each rank $k \geq 3$, there is a string C-group representation of rank $k$ with group $A_n$ for some $n$.  In particular, for each even rank $r \geq 4$, there is a string C-group representation of Schl\"afli type $\{10,3^{r-2}\}$ with group $A_{2r+1}$, and for each odd rank $q \geq 5$, there is a string C-group representation of Schl\"afli type $\{10,3^{q-4},6,4\}$ with group $A_{2q+3}$.
\end{theorem}

\begin{table}
\begin{small}
\begin{center}
\begin{tabular}{||c|c|c||}
\hline
Group&Schl\"afli Type&CPR Graph\\
\hline%evenrank
\begin{tabular}{c}
\\[-10pt]
$A_{2r+1}$\\[5pt] $^{(r\textrm{ even and }\geq 4)}$
\end{tabular} & $\{10,3^{r-2}\}$
&
$\xymatrix@-1.7pc{&& *+[o][F]{}  \ar@{-}[rr]^1 && *+[o][F]{}  \ar@{-}[rr]^2 && *+[o][F]{}  \ar@{-}[rr]^3 && *+[o][F]{}   \ar@{.}[rr] && *+[o][F]{}  \ar@{-}[rr]^{r-2} && *+[o][F]{}  \ar@{-}[rr]^{r-1} && *+[o][F]{} \\
&& && && && && && && \\
 *+[o][F]{}  \ar@{-}[rr]_0 && *+[o][F]{}  \ar@{-}[rr]_1 && *+[o][F]{}   \ar@{-}[uu]^0   \ar@{-}[rr]_2 && *+[o][F]{}  \ar@{-}[uu]^0   \ar@{-}[rr]_3 && *+[o][F]{}  \ar@{-}[uu]^0     \ar@{.}[rr] && *+[o][F]{}  \ar@{-}[uu]_0   \ar@{-}[rr]_{r-2} && *+[o][F]{}  \ar@{-}[uu]_0   \ar@{-}[rr]_{r-1} &&  *+[o][F]{} \ar@{-}[uu]_0     }$\\
\hline%oddrank
\begin{tabular}{c}
\\[-10pt]
$A_{2r+3}$\\[5pt] $^{(r \textrm{ odd and } \geq 5)}$
\end{tabular} & $\{10,3^{r-4},6,4\}$
&
$ \xymatrix@-1.7pc{&& *+[o][F]{}  \ar@{-}[rr]^1 && *+[o][F]{}  \ar@{-}[rr]^2 && *+[o][F]{}  \ar@{-}[rr]^3 && *+[o][F]{}   \ar@{.}[rr] && *+[o][F]{}  \ar@{-}[rr]^{r-2} && *+[o][F]{}  \ar@{-}[rr]^{r-1} && *+[o][F]{}   \ar@{-}[rr]^{r-2} && *+[o][F]{}\\
&& && && && && && && &&\\
 *+[o][F]{}  \ar@{-}[rr]_0 && *+[o][F]{}  \ar@{-}[rr]_1 && *+[o][F]{}   \ar@{-}[uu]^0   \ar@{-}[rr]_2 && *+[o][F]{}  \ar@{-}[uu]^0   \ar@{-}[rr]_3 && *+[o][F]{}  \ar@{-}[uu]^0     \ar@{.}[rr] && *+[o][F]{}  \ar@{-}[uu]_0   \ar@{-}[rr]_{r-2} && *+[o][F]{}  \ar@{-}[uu]_0   \ar@{-}[rr]_{r-1} &&  *+[o][F]{} \ar@{-}[uu]_0  \ar@{-}[rr]_{r-2} && *+[o][F]{}   \ar@{-}[uu]_0  }$\\
 \hline
\end{tabular}
\caption{String C-group representations for $A_n$}\label{TT}
\end{center}
\end{small}
\end{table}

We then managed to construct string C-group representations of rank $\lfloor \frac{n-1}{2} \rfloor$ for $A_n$ with $n\geq 12$.

\begin{theorem}\cite[Theorem 1.1]{flm2}
For each $n \notin \{3,4,5,6,7,8,11\}$, there is a rank $\lfloor \frac{n-1}{2} \rfloor$ string C-group representation of the alternating group $A_n$.
\end{theorem}

We were quickly convinced that this rank was the best possible but it took another five years, and Peter Cameron joining forces, to finally prove the following theorem.

\begin{theorem}\cite[Theorem 1.1]{cflm}\label{cflmTheo}
The maximum rank of a string C-group representation of $A_n$ is $3$ if $n=5$; $4$ if $n=9$; $5$ if $n=10$; $6$ if $n=11$ and $\lfloor\frac{n-1}{2}\rfloor$ if $n\geq 12$.
Moreover, if $n=3, 4, 6, 7$ or $8$, the group $A_n$ is not a string C-group. 
\end{theorem}

Fernandes and the author then managed to construct string C-group representations of each rank $3\leq r \leq \lfloor (n-1)/2 \rfloor$ for $A_n$ with $n\geq 12$.

\begin{theorem}\cite[Theorem 1.1]{fl2}\label{rralt}
For $n\geq 12$ and  for every $3 \leq r \leq \lfloor (n-1)/2 \rfloor$, the group $A_n$ has at least one string C-group representation of rank $r$.
\end{theorem}

Table~\ref{an} together with the proof of Theorem~\ref{cflmTheo} suggests that the enumeration results we obtained for the symmetric groups will be much harder to get for the alternating groups. Nevertheless, this is a very interesting problem as well so we list it here.

\begin{problem}\label{anenum}
Enumerate the string C-group representations of rank $r$ of $A_n$ for every $3\leq r \leq n-1$.
\end{problem}

As for the symmetric groups, a first step in trying to solve Problem~\ref{anenum} has been made by Kiefer and the author. We 
computed the number of unordered pairs of commuting involutions in $A_n$, up to conjugacy in $Aut(S_n)$ in~\cite{KL2013}. We obtained the following result.

\begin{theorem}\cite[Theorem 1.2]{KL2013}\label{pairalt}
Let $n>1$ be a positive integer.
Define $\lambda(k)$, $\phi(k,n)$ and $\mu(n)$ as follows.
\begin{equation*}\begin{aligned}
\lambda(k) =\Big\lfloor{\Big(\frac{k}{2}+1\Big)^2}\Big\rfloor
\end{aligned}\end{equation*}
\begin{equation*}
\begin{aligned}
\phi(k,n) = \left\{ \begin{aligned}
& \lambda(k) - 1 & \textrm{ if } k \leq \Big\lfloor\frac{n}{2}\Big\rfloor,\\
& \lambda(k) - \psi(k,n) -1 & \textrm{ if } k > \Big\lfloor\frac{n}{2}\Big\rfloor,
\end{aligned}\right.
\end{aligned}\end{equation*}

\begin{equation*}\begin{aligned}
\mu(n) = & -2\Big\lfloor\frac{n}{2}\Big\rfloor+\sum_{\substack{k=1 \\ k \textrm{ even}}}^{n} \left[ 
\gamma(k)\cdot \Big\lceil \frac{1}{2} \cdot \Big(n -k+ 1\Big) \Big\rceil + 
\delta(k)\cdot \Big\lfloor\frac{1}{2} \cdot \Big(n-k+ 1\Big)\Big\rfloor \right]
\end{aligned}\end{equation*} where
\begin{equation*}\begin{aligned}
\gamma(k) & = \frac{k^2}{8} + \frac{3k}{4} + 1,\\
\delta(k) & =  \frac{k^2}{8} + \frac{k}{4}.
\end{aligned}\end{equation*}
There are, up to isomorphism,
\begin{equation*}\begin{aligned}
\frac{1}{2} \Big(\mu(n) + \sum_{\substack{k=1 \\ k \textrm{ even}}}^{n} \phi(k,n) \Big)
\end{aligned}\end{equation*}

unordered pairs of commuting involutions in $A_{2n}$ and $A_{2n+1}$ except for $A_6$ in which there is, up to isomorphism, a unique 
unordered pair of commuting involutions.
\end{theorem}
As mentioned in the previous section, Table~\ref{snancommut} gives the number of pairwise nonisomorphic pairs of commuting involutions for $A_n$ for some values of $n$.

\section{Projective linear groups}
Several projective linear groups were analyzed in~\cite{LV}. The observation of the data collected permitted, over the years, to obtain classification results that we summarize in this section.
\subsection{Groups $PSL(2,q)$} The values of $q$ for which
a $PSL(2,q)$ group has rank three polytopes were determined by Sjerve and Cherkassoff.
\begin{theorem}\cite[Theorem 1.3]{SC94}
The $PSL(2,q)$ group may be generated by three involutions, two of which commute, if and only if $q\neq 2, 3, 7$ or $9$.
\end{theorem}

For the groups $PSL(2,q)$, we quickly observed that the maximal rank of a string C-group representation for these groups is 4 as they do not possess subgroups that are direct products of two dihedral groups of order at least 6 each, obtaining the following theorem.

\begin{theorem}~\cite[Theorem 2]{LV}
Let $G \cong PSL(2,q)$. The rank of a string C-group representation of $G$ is at most 4.
\end{theorem}

More striking was the fact that, up to $q=32$, only two groups $PSL(2,q)$ had a rank four representation, namely $PSL(2,11)$ and $PSL(2,19)$.
Together with Schulte, we managed to prove the following theorem.

\begin{theorem}\label{ourthm}\cite[Theorem 1]{ls07}
If $PSL(2,q)$ has a string C-group representation of rank $4$, then $q=11$ or $19$. Moreover, for each of these two values of $q$, there is, up to isomorphism, a unique string C-group representation of rank four for $PSL(2,q)$.
\end{theorem}
This result made the 11-cell and the 57-cell even more special. They were the only two string C-group representations of rank four known for groups $PSL(2,q)$ and our result explained why.

The proof of Theorem~\ref{ourthm} and subsequent theorems on groups with socle $PSL(2,q)$ (see next subsections) rely heavily on the complete knowledge of the subgroup structure of $PSL(2,q)$.
As mentioned in~\cite[Section 4.1]{BDL2010}, the subgroup structure of $PSL(2, q)$ was
first obtained in papers by Moore~\cite{Moo1904} and Wiman~\cite{Wiman1899}.

Observe also that the more recent results on $PSL(4,q)$ groups make it less surprising that the maximal rank of a string C-group representation of $PSL(2,q)$ is three in most cases. Indeed, there is not much freedom left by the commuting property to find larger sets of involutions that would give string C-group representations for these groups.
\subsection{Groups $PGL(2,q)$}
Again, in the rank three case, Sjerve and Cherkassoff  determined for which values of $q$ the group $PGL(2,q)$ has a string C-group representation.
\begin{theorem}\cite[Theorem 1.4]{SC94}
The $PGL(2,q)$ group may be generated by three involutions, two of which commute, if and only if $q\neq 2$.
\end{theorem}
Looking at the data collected in~\cite{LV}, we conjectured with Schulte that the maximum rank of a string C-group representation of $PGL(2,q)$ should be three except when $q=5$ and we proved the following result.
\begin{theorem}\label{ourthm2}\cite[Theorem 4.8]{ls09}
The group $PGL(2,q)$ has string C-group representations of rank at most 4. Moreover, only $PGL(2,5)$ has a string C-group representation of rank 4.
\end{theorem}
The only exception in rank four is the 4-simplex whose automorphism group is $S_5\cong PGL(2,5)$.
\subsection{Almost simple groups with socle $PSL(2,q)$}
Later on, with Connor and De Saedeleer, we decided to study almost simple groups with socle $PSL(2,q)$. We proved the following classification theorem, in the same vein as the results obtained with Schulte, answering at the same time a conjecture we stated in~\cite{ls09}.

\begin{theorem}
\label{dijutho}\cite[Theorem 1.1]{DiJuTho}
Let $\PSL(2,q)\leq G \leq \PGaL(2,q)$. Suppose $G$ has a string C-group representation.
Then 
\begin{enumerate}
 \item if $q=2$ then $G \cong \PSL(2,2) \cong S_3$ and $G$ has a unique rank $2$ string C-group representation, namely the one coming from the triangle;
 \item if $q=3$ then $G\cong \PGL(2,3)\cong\Sym_4$ and $G$ only has rank three string C-group representations;
 \item if $q=4$ or $5$ then either $G\cong\PSL(2,4)\cong \PSL(2,5)\cong\Alt_5$ and $G$ has  rank three string C-group representations only, or $G\cong \PGL(2,5)\cong\Sym_5$ and $G$ has rank three and four string C-group representations;
 \item if $q=7$ then $G\cong \PGL(2,7)$ and $G$ has rank three string C-group representations only;
 \item if $q\geq 8$ then
 \begin{enumerate}
  \item if $q=2^{2k+1}$, $k\geq 1$, then $G\cong \PSL(2,2^{2k+1})$ and $G$ has rank three string C-group representations only;
 \item if $q=9$ then either $G\cong \PGL(2,9)$, or $G\cong \PSiL(2,9)\cong \Sym_6$, or $G\cong\PGaL(2,9)$, and $G$ has rank three string C-group representations;
 moreover $\PSiL(2,9)$ has string C-group representations of ranks $3$, $4$ and $5$;
  \item if $q=p^{2k+1}\geq11$, $p$ an odd prime and $k\geq 0$, then $G\cong \PSL(2,p^{2k+1})$ or $G\cong \PGL(2,p^{2k+1})$; in either case, $G$ has rank three string C-group representations;
  if moreover $q=11$ or $19$ then $G\cong \PSL(2,q)$ has rank four string C-group representations;
  \item if $q=p^{2k}\geq 16$, $p$ any prime and $k\geq 1$ then either $G\cong\PSL(2,p^{2k})$ or $G\cong\PGL(2,p^{2k})$ or  $G\cong \PSL(2,p^{2k})\rtimes\langle\beta\rangle$ or $G\cong\PGL(2,p^{2k})\rtimes\langle\beta\rangle$, where $\beta$ is a Baer involution of $\PGaL(2,p^{2k})$; in all four cases, $G$ has rank three string C-group representations;
   moreover, $\PSL(2,p^{2k})\rtimes\langle\beta\rangle$ also has string C-group representations of rank $4$.
 \end{enumerate}
\end{enumerate} 
\end{theorem}
\subsection{Open problems on groups with socle $PSL(2,q)$}

Enumeration of C-group representations of rank three for groups $PSL(2,q)$ and $PGL(2,q)$ can be found in~\cite{CPS2008}. Even though their results are for hypermaps, all the hypermaps they give satisfy the intersection property by~\cite[Corollary 4.2]{CO2013} as none of $PSL(2,q)$ or $PGL(2,q)$ possesses a non-trivial cyclic subgroup that is normal. 
%Their results being quite long to summarize, we refer the reader to their paper.
String C-group representations of rank four for $PSL(2,q)$ and $PGL(2,q)$ are known thanks to Theorems~\ref{ourthm} and~\ref{ourthm2}.
\begin{problem}
Enumerate the non-isomorphic string C-group representations of rank three and four of $\PSL(2,p^{2k})\rtimes\langle\beta\rangle$ where $\beta$ is a Baer involution of $\PGaL(2,p^{2k})$.
\end{problem}
\begin{problem}
Enumerate the non-isomorphic string C-group representations of rank three of 
$\PGL(2,p^{2k})\rtimes\langle\beta\rangle$ where $\beta$ is a Baer involution of $\PGaL(2,p^{2k})$.
\end{problem}
\subsection{Groups $PSL(3,q)$ and $PGL(3,q)$.}
In~\cite{BV2010}, Brooksbank and Vicinsky proved the following theorem.

\begin{theorem}\cite{BV2010}
If $G\leq GL(3,q)$ has a string C-group representation, then $q$ is odd and there is a non-degenerate symmetric bilinear form $\bf f$ on $V$ such that $\Omega(V,{\bf f})\leq G \leq I(V,{\bf f})$.
\end{theorem}
As a corollary of their result (see~\cite{Adrien}), the groups $PSL(3,q)$ and $PGL(3,q)$ are not string C-groups.

The proof of the Brooksbank-Vicinsky theorem is very different from the ones for the groups with socle $PSL(2,q)$. It uses the fact that, for any subgroup $G$ of $GL(3,q)$, any string C-group representation of $G$ is of rank at most four. Then, the authors prove that if $G$ has a string C-group representation of rank three or four, $q$ is odd and $G$ preserves a non-degenerate symmetric bilinear form. 
%This result was extended later in dimension 4 (see Section~\ref{section:psl4}).
\subsection{Open problems on groups with socle $PSL(3,q)$}
From the data collected in~\cite{LV}, it is clear that some of these groups have string C-group representations of rank three, four and five. Again, it would be nice to know what is the maximal rank for these groups and to enumerate them.

\begin{problem}
Determine the maximal rank of a string C-group representation for $G$ with $PSL(3,q)\leq G \leq Aut(PSL(3,q))$.
\end{problem}
\begin{problem}
Enumerate the non-isomorphic string C-group representations of rank three of $G$ with $PSL(3,q)\leq G \leq Aut(PSL(3,q))$.
\end{problem}
\subsection{Groups with socle $PSL(4,q)$}\label{section:psl4}
In~\cite{BL2015}, with Brooksbank, we showed that the groups $PSL(4,q)$ have string C-group representations if and only if $q$ is odd. Moreover, we showed that for each odd $q$, there is at least one string C-group representation of rank four.

\begin{theorem}
\label{psl4-odd}\cite[Theorem 1.1 and Corollary 4.5]{BL2015}
If $q=p^k$ for odd $p$, then
$PSL(4,q)$ has at least one string C-group representation of rank $4$.
If $q$ is even, then $PSL(4,q)$ has no string C-group representation.
\end{theorem}

\begin{conjecture}
The group $PSL(4,q)$, with $q$ odd, has maximal rank 4 as a string C-group.
\end{conjecture}

\begin{problem}
Determine the maximal rank of a string C-group representation for $G$ with $PSL(4,q)\leq G \leq Aut(PSL(4,q))$.
\end{problem}

\section{Suzuki groups}\label{suzuki}
Among the nonabelian simple groups, those with the easiest subgroup structure are the Suzuki simple groups. They are therefore the most promising groups for which the three main questions mentioned in the introduction seem solvable. This is the reason why we decided to start working on these questions with these groups. As a rule of thumb, it seems to us that if a general question on finite simple groups is not answerable for Suzuki simple groups, it is unlikely it will be for the other families.

We recall the definition of the Suzuki groups as given in \cite{Lun80}.
Let $\mathcal K$ be a field of characteristic 2 with $\mid {\mathcal K} \mid > 2$.
Let $\sigma$ be an automorphism of $\mathcal K$ such that $x^{\sigma^2} = x^2$ for each $x$ in $\mathcal K$.
Let $\mathcal B$ be the 3-dimensional projective space over $\mathcal K$ and let $(x_0, x_1, x_2, x_3)$ be the coordinates of a point of $\mathcal B$.
Let $E$ be the plane defined by the equation $x_0$ = 0 and put $U = (0, 1, 0, 0) \mathcal K$.
We introduce coordinates in the affine space ${\mathcal B}_E$ by $x = \frac{x_2}{x_0}$,
$y = \frac{x_3}{x_0}$, and $z = \frac{x_1}{x_0}$.
Finally, let $\mathcal D$ be the set of points of $\mathcal B$ consisting of $U$ and all those points of ${\mathcal B}_E$
whose coordinates $(x, y, z)$ satisfy the equation
\[z = xy + x^{\sigma + 2} + y^\sigma \]
We denote by {\it Sz($\mathcal K$, $\sigma$)} the group of all projective collineations of $\mathcal B$ which leave $\mathcal D$ invariant.
When $\mathcal K$ is a finite field, $\sigma$ is unique and in this case, $\mathcal K$ is isomorphic to {\it GF(q)} with $q = 2^{2e+1}$.
The groups {\it Sz(q)} are the {\it Suzuki groups} named after Michio Suzuki who found
them in 1960. The generalizations {\it Sz($\mathcal K$, $\sigma$)} are due to Rimhak Ree and Jacques Tits (see for example \cite{Tit62b}).
The set $\mathcal D$ is an {\it ovoid}, i.e. a non-empty point-set of a projective 3-space that satisfies the following three conditions.
\begin{enumerate}
\item No three points are collinear.
\item If $p \in\mathcal D$, there exists a plane $E$ of $\mathcal B$ with ${\mathcal D} \cap E = \{ p \}$.
\item If $p \in\mathcal D$ and if $E$ is a plane of $\mathcal B$ with ${\mathcal D}\cap E = \{p\}$, then all lines
{\it l} through $p$ which are not contained in $E$ carry a point of $\mathcal D$ distinct from $p$.
\end{enumerate}

Every involution of $Sz(q)$ fixes a unique point of $\mathcal D$. Moreover, commuting involutions have the same fixed point.
Looking at the results obtained in~\cite{LV}, and thanks to the knowledge gathered on Suzuki simple groups during our PhD thesis, we quickly managed to prove the following theorem.

\begin{theorem}\label{sz}~\cite[Theorem 1]{Lee2005a}
Let $Sz(q) \leq G \leq Aut(Sz(q))$ with $q = 2^{2e+1}$ and $e > 0$ a positive integer. Then $G$ is a C-group if and only if $G = Sz(q)$. Moreover, if $(G, \{\rho_0,\ldots, \rho_{n-1}\})$ is a string C-group, then $n=3$.
\end{theorem}
We give here a different proof from the original one, using centralisers of involutions.
\begin{proof}
The group $G:=Sz(q)$ has a unique conjugacy class of involutions. Suppose $G$ is the natural permutation representation of $Sz(q)$ acting on a Suzuki-Tits ovoid.
Two involutions commute in $G$ if and only if they fix the same point on the ovoid. As $G$ is simple, it implies that the maximal rank of a string C-group representation of $G$ is three. Moreover, since $Aut(Sz(q)) \cong Sz(q):C_{2e+1}$ and $2e+1$ is odd, all involutions of $Aut(Sz(q))$ are in $Sz(q)$ and therefore, a set of involutions of $Aut(Sz(q))$ can at most generate $Sz(q)$.
Finally, take two involutions $\rho_0$ and $\rho_2$ that commute in $G$. Pick an involution $\rho_1$ such that the order of $\rho_0\rho_1$ is $q-1$. Then $\rho_1$ obviously does not commute with $\rho_2$ and $G = \langle \rho_0,\rho_1,\rho_2\rangle$ as $\langle \rho_0,\rho_1\rangle$ is a maximal subgroup of $G$. The pair $(G,\{ \rho_0,\rho_1,\rho_2\})$ is a string C-group representation of $G$.
\end{proof}

In 2010, with Ann Kiefer, we managed to count the number of pairwise nonisomorphic string C-group representations of rank three of a given Suzuki group $Sz(q)$.

\begin{theorem}\label{polys}~\cite[Theorem 2]{KL2010}
Up to isomorphism and duality, a given Suzuki group $\mathrm{Sz}(q)$, with $q=2^{2e+1}$ and $e >0$ an integer, has \[\frac{1}{2} \sum_{2f+1 \vert 2e+1} \mu(\frac{2e+1}{2f+1}) \sum_{\substack{n \vert 2f+1 \\ n \neq 1}} \lambda(n) \psi(n,2f+1)\] string C-group representations, where $\mu(n)$ is the Moebius function,
\begin{align*}
\lambda(n) & = \frac{1}{n} \sum_{d \mid n} \mu(\frac{n}{d}) \cdot 2^d \text{ and} \\
\psi(n,2f+1)&  = \sum_{m \mid \frac{2f+1}{n}}\frac{ \sum_{d \mid m} \mu(\frac{m}{d})( 2^{nd}-1)}{m}.
\end{align*}
All these representations are non-degenerate, i.e. have a Schl\"afli type with entries $\geq 3$.
\end{theorem}
This result closed the chapter of string C-group representations for Suzuki simple groups as the three main questions we gave in Section~\ref{intro} were answered.

\section{Small Ree groups}

The small Ree groups $\Ree(q)$, defined over a finite field of order $q=3^{2e+1}$ and $e>0$, were discovered by Rimhak Ree~\cite{Ree1960} in 1960. 
In the literature they are also denoted by $\Reeb(q)$. These groups have a subgroup structure quite similar to that of the Suzuki simple groups $\Sz(q)$, with $q=2^{2e+1}$ and $e>0$.  Suzuki and Ree groups play a somewhat special role in the theory of finite simple groups, since they exist because of a Frobenius twist arising from a special automorphism of the field over which they are constructed.
They have no counterpart in characteristic zero. 

The Ree group $G:=\Ree(q)$, with $q=3^{2e+1}$ and $e\geq 0$, is a group of order $q^3(q-1)(q^3+1)$. It has a natural permutation representation on a Steiner system ${\mathcal S} := (\Omega,\mathcal{B}) = S(2,q+1,q^3+1)$ consisting of a set $\Omega$ of $q^3+1$ elements, the {\em points\/}, and a family of $(q+1)$-subsets $\mathcal{B}$ of $\Omega$, the {\em blocks\/}, such that any two points of $\Omega$ lie in exactly one block. This Steiner system is also called a {\em Ree unital\/}. In particular, $G$ acts $2$-transitively on the points and transitively on the incident pairs of points and blocks of ${\mathcal S}$.  

A list of the maximal subgroups of $G$ is available, for instance, in~\cite[p. 349]{VM2000} and~\cite{Kle1988}.

The group $G$ has a unique conjugacy class of involutions (see~\cite{Ree1960}). Every involution $\rho$ of $G$ has a  block $B$ of $\mathcal S$ as its set of fixed points, and $B$ is invariant under the centralizer $C_{G}(\rho)$ of $\rho$ in $G$. Moreover, $C_G(\rho)\cong C_2\times PSL(2,q)$, where $C_2 = \langle \rho \rangle$ and the $PSL(2,q)$-factor acts on the $q+1$ points in $B$ as it does on the points of the projective line $PG(1,q)$. Hence the knowledge of string C-group representations of groups $PSL(2,q)$ is helpful to construct string C-group representations of large rank for $G$.

The automorphism group $\mathsf{Aut}(\Ree(q))$ of $\Ree(q)$ is given by
\[ \mathsf{Aut}(\Ree(q)) \cong \Ree(q)\!:\!C_{2e+1}, \]
so in particular $\mathsf{Aut}(\Ree(3)) \cong \Ree(3)$. 

In~\cite{Leemans:2015}, using the list of maximal subgroups of $G$ and the geometric properties of the Steiner system, we obtained, with Egon Schulte and Hendrik Van Maldeghem, the following theorem bounding the rank of a string C-group representation of a small Ree group $\Ree(q)$.

\begin{theorem}\label{smallree}\cite[Theorem 1.1]{Leemans:2015}
Among the almost simple groups $G$ with $\Ree(q) \leq G \leq \mathsf{Aut}(\Ree(q))$ and $q = 3^{2e+1}\neq 3$, only the Ree group $\Ree(q)$ itself is a C-group. In particular, $\Ree(q)$ admits a representation as a string C-group of rank $3$, but not of higher rank. Moreover, the non-simple Ree group $\Ree(3)$ is not a C-group.
\end{theorem}

In other words, the groups $\Ree(q)$ behave just like the Suzuki groups:\ they allow representations as string C-groups, but only of rank $3$. The proof of this theorem, unlike its counterpart for the Suzuki groups, is very lengthy and uses a deep analysis of the lattice of subgroups of a Ree group.

\begin{problem}
Enumerate the non-isomorphic string C-group representations of rank three of $\Ree(q)$.
\end{problem}

\section{Orthogonal and symplectic groups}

With Brooksbank and Ferrara, we decided to investigate the orthogonal groups and to look for string C-group representations of large rank. Indeed, apart from the symmetric and alternating groups, as well as some crystallographic groups studied by Barry Monson and Egon Schulte (see~\cite{MS1, MS2, MS3}), no other family of groups was known to have possible large ranks that could grow linearly depending on one of its parameters, being the permutation degree or the dimension of the group or the size of the field on which the group is defined.
We obtained the following theorem for orthogonal and symplectic groups.

\begin{theorem}\cite[Corollary 1.3]{BFL2018}
\label{coro:main-sp}
For each integer $k\geq 2$, positive integer $m$, and $\epsilon\in\{-,+\}$,
the orthogonal group $\Orth^{\epsilon}(2m,\Fk)$ has a string C-group representation of rank $2m$, and the symplectic group $\Sp(2m,\Fk)$ has a string C-group representation of rank $2m+1$.
\end{theorem}

Given the way we had to choose the involutions to construct the string C-group representations mentioned in Theorem~\ref{coro:main-sp}, we have good reasons to believe that no higher rank can be achieved for these groups.

\begin{problem}
Determine the maximal rank of a string C-group representation for $\Orth^{\epsilon}(2m,\Fk)$ and for $\Sp(2m,\Fk)$.
\end{problem}

\section{Sporadic groups}
In~\cite{LV}, all string C-group representations were computed for the sporadic groups
$M_{11}$,
$M_{12}$,
$M_{22}$,
$J_{1}$
and $J_{2}$, as well as for their respective automorphism groups.
In 2010, Hartley and Alexander Hulpke~\cite{HHalg} designed more efficient algorithms that permitted them to classify all string C-group representations of the five Mathieu groups, the first three Janko groups, the Higman-Sims group, the McLaughlin group and the Held group.
In 2012, with Mixer, we further improved the algorithms and managed to classify all string C-group representations of the third Conway group~\cite{LM2012}.
In~\cite{CLM2012}, we proved with Connor and Mixer that the maximal rank of a string C-group representation for the O'Nan group is four. Moreover, we gave all string C-group representations of rank four for the O'Nan group.
In~\cite{CL2016}, with Connor, we enumerated the regular maps of the O'Nan group using character theory, but we were unable to determine how many of them give a string C-group representation of rank three of the O'Nan group.
Finally, in~\cite{LM2020}, with Jessica Mulpas, we classified all string C-group representations of the Rudvalis group and the Suzuki group.

A summary of these results is given in Table~\ref{sporadics}. No string C-group representation of rank six or higher exists for all the groups listed in that table as shown in the corresponding references given above.

\begin{table}
\begin{tabular}{|c|c|c|c|c|}
\hline
$G$&Order of $G$&Rank 3&Rank 4&Rank 5\\
\hline
$M_{11}$&7,920&0&0&0\\
$M_{12}$&95,040&23&14&0\\
$M_{22}$&443,510&0&0&0\\
$M_{23}$&10,200,960&0&0&0\\
$M_{24}$&244823040&490&155&2\\
\hline
$J_1$&175,560&148&2&0\\
$J_2$&604,800	&137&17&0\\
$J_3$&50,232,960&303&2&0\\
\hline
$HS$&44,352,000&252&57&2\\
$McL$& 898,128,000&0&0&0\\
$He$ &4,030,387,200 &1188&76& 0\\
$Ru$&145,926,144,000&21594&227&0\\
$Suz$&448,345,497,600&7119&257&13\\
$O'N$&460,815,505,920&Unknown&16&0\\
$Co_3$&495,766,656,000&10586&873&22\\
\hline
\end{tabular}
\caption{Number of string C-group representations for sporadic groups}\label{sporadics}
\end{table}

\begin{problem}
Determine the number of pairwise nonisomorphic string C-group representations of rank three for the O'Nan group.
\end{problem}

\begin{problem}
Try to push further the algorithms to study the remaining sporadic groups.
\end{problem}
One of the problems of the current algorithms is that they work with a permutation representation of the group to be analyzed. This becomes a problem when such a representation is on more than, say, 150000 points. Using matrix groups instead of permutation groups might help.
Another way to improve the existing algorithms would be to use parallel computing.

\section{Collateral results}

While working on specific families of simple groups, more general results were found. We mention the main ones in this section.

\subsection{Transitive groups}

Transitive groups played a key role in the proof of Theorem~\ref{cflmTheo}.
Indeed, to prove that theorem, the authors had to get the best possible bound for the maximal rank of a string C-group representation for a transitive group of degree $n$ that is not $S_n$ nor $A_n$. Together with Cameron, Fernandes and Mixer, we obtained the following result.

\begin{theorem}\cite[Theorem 1.2]{cflm2}
Let $\Gamma$ be a string C-group representation of rank $d$ which is isomorphic to a
transitive subgroup of $\Sym_n$ other than $\Sym_n$ or $\Alt_n$. Then one of
the following holds:
\begin{enumerate}
\item $d\le n/2$;
\item $n \equiv 2 \mod{4}$, $d = n/2+1$ and $\Gamma$ is $\Cyc_2\wr \Sym_{n/2}$. The generators are
\begin{center}
$\rho_0 = (1,n/2+1)(2,n/2+2)\ldots (n/2,n)$;

$\rho_1 = (2,n/2+2)\ldots (n/2,n)$;

$\rho_i = (i-1,i)(n/2+i-1,n/2+i)$ for $2\leq i \leq n/2$.
\end{center}
Moreover the Schl\"afli type is $\{2,3, \ldots, 3,4\}$. 
\item $\Gamma$ is transitive imprimitive and is one of the examples appearing in Table~\ref{ploys}.
\item $\Gamma$ is primitive. In this case, $n=6$. Moreover $\Gamma$ is obtained from the permutation representation of degree $6$ of $\Sym_5 \cong \PGL_2(5)$ and it is the group of the $4$-simplex of Schl\"afli type $\{3,3,3\}$.

\begin{table}
\begin{center}
\begin{tabular}{|c|c|c|c|l|}
\hline
$Degree$&$Number$&Structure&Order&Schl\"afli type\\
\hline
\hline
6&9&$\Sym_3 \times \Sym_3$&36&$\{3,2,3\}$\\
\hline
6&11&$2^3:\Sym_{3}$&48&$\{2,3,3\}$\\
6&11&$2^3:\Sym_{3}$&48&$\{2,3,4\}$\\
\hline
8&45&$2^4:\Sym_{3}:\Sym_{3}$&576&$\{3,4,4,3\}$\\
\hline
\end{tabular}
\caption{Examples of transitive imprimitive string C-groups of degree $n$ and rank $n/2+1$ for $n\leq 9$.}\label{ploys}
\end{center}
\end{table}

\end{enumerate}
%In case (a), if equality holds, then there is a unique example for each $n\equiv2$ (mod~$4$).
\end{theorem}

This result and the experimental data obtained for the symmetric groups make us strongly believe in Conjecture~\ref{conjsn}. Indeed, it forces all the maximal parabolic subgroups of a string C-group representation of $S_n$ of rank at least $n/2+3$ to be intransitive subgroups of $S_n$ when $n$ is large enough, that is $n\geq 9$.

\subsection{Rank reduction}

The proof of Theorem~\ref{rralt} has  inspired Peter Brooksbank and the author to prove a rank reduction theorem~\cite{BL2018}.

\begin{theorem}\cite[Theorem 1.1]{BL2018}\label{thm:main}
Let $(G;\{\rho_0,\ldots, \rho_{n-1}\})$ be an irreducible string C-group of rank $n\geq 4$.
If $\rho_0\in \langle \rho_0\rho_2,\rho_3\rangle$, then 
$(G;\{\rho_1, \rho_0\rho_2, \rho_3, \ldots, \rho_{n-1}\})$ is a string C-group of rank $n-1$.
\end{theorem}

In particular, this theorem gives the following corollary.

\begin{corollary}\cite[Corollary 1.3]{BL2018}
\label{coro:odd-order-sequence}
Let $(G;\{\rho_0,\ldots, \rho_{n-1}\})$ be an irreducible string C-group representation of rank $n\geq 4$ of $G$.
Let $\{p_1,\ldots, p_{n-1}\}$ be its Schl\"afli type, and put
\[
t=\max\{j\in\{0,\ldots,n-3\}\colon \forall i\in\{0,\ldots,j\},\; p_{2+i}\;\mbox{is odd}\}.
\]
Then $G$ has a string C-group representation of rank $n-i$ for each $i\in\{0,\ldots,t\}$.
\end{corollary}

This corollary stresses the importance of trying to construct string C-group representations of large ranks as they may give representations of lower ranks for free.
For instance, the corollary makes the proof of Theorem~\ref{1} straightforward. 
It can also be used to improve the proof of Theorem~\ref{rralt}, and it was used by Brooksbank and the author to prove the following result.

\begin{theorem}
\cite[Theorem 1.4]{BL2018}\label{thm:sp}
Let $k\geq 2$ and $m\geq 2$ be integers.
\begin{enumerate}
\item[(a)]  The symplectic group $\Sp(2m,\Fk)$ has a string C-group representation of rank $n$ for each $3\leq n\leq 2m+1$.
 \item[(b)]  The orthogonal groups $\Orth^{+}(2m,\Fk)$ and $\Orth^-(2m,\Fk)$ have
 string C-group representations of rank $n$ for each $3\leq n\leq 2m$.
 \end{enumerate}
\end{theorem}

We conclude this section by recalling the following conjecture made by Brooksbank and the author in~\cite{BL2018}.

\begin{conjecture}\cite[Conjecture 5.1]{BL2018}
The group $A_{11}$ is the only finite simple group whose set of ranks of string C-group 
representations is not an interval in the set of integers.
\end{conjecture}

\section{C-groups}
A natural question is to see if one could prove similar results by weakening the hypotheses. We finish this survey with what is currently known about C-group representations of almost simple groups.

Work with Thomas Connor showed that the word `string' in the last part of the statement of the Theorem~\ref{sz} can be removed.

\begin{theorem}\label{cgroupssz}\cite[Classification Theorem 1.1]{CL2015}
Let $q=2^{2e+1}\neq 2$ be an odd power of $2$ and let $Sz(q)\leq G \leq Aut(Sz(q))$ be an almost simple group of Suzuki type. Let $(G,\{\rho_0,\ldots,\rho_{r-1}\})$ be a $C$-group representation of $G$. Then $G=Sz(q)$ and $r=3$. Moreover there exist at least one string C-group representation and one nonstring C-group representation of $G$.
\end{theorem}

C-group representations in general are also very interesting to study as they are smooth quotients of Coxeter groups, not only those with a string diagram. Moreover, they are often (but not always as in the case of string C-group representations) automorphism groups of more general geometric objects that, together with Fernandes and Weiss, we called hypertopes in~\cite{flw}. These objects are apartments in the theory of Buildings due to Jacques Tits (see~\cite{Tit74}), one more reason to study them further.

It would be interesting to have a result similar to Theorem~\ref{polys} for C-groups.
\begin{problem}
Determine the number of pairwise non-isomorphic C-group representations of $Sz(q)$.
\end{problem}

%\subsection{C-group representations of $PSL(2,q)$ and $PGL(2,q)$}
Together with Connor and Sebastian Jambor, we obtained the following result for groups $PSL(2,q)$ and $PGL(2,q)$.

\begin{theorem}\label{theorem:main}~\cite[Theorem 1.1]{CJL2015}
Let $G\cong \PSL(2,q)$ for some prime power $q \geq 4$. A C-group representation of $G$ is of rank $4$ if and only if $q\in\{7,9,11,19,31\}$. Otherwise it is $3$. 

Let $G\cong \PGL(2,q)$ for some prime power $q \geq 4$. A C-group representation of $G$ is of rank $4$ if and only if $q=5$. Otherwise it is $3$. 
\end{theorem}
Our proof was in two steps: first we reduced the possible Coxeter diagrams of the C-group representations of group $PSL(2,q)$ and $PGL(2,q)$ using the subgroup structure of $PGL(2,q)$. Then we used the $L_2$-quotient algorithm to see which diagrams gave Coxeter groups admitting groups $PSL(2,q)$ or $PGL(2,q)$ as quotients. Unfortunately, this algorithm does not permit to recognize other almost simple groups with socle $PSL(2,q)$ so we were not able to extend our result to all almost simple groups with socle $PSL(2,q)$.

\begin{problem}
Prove a theorem similar to Theorem~\ref{dijutho} for C-groups.
\end{problem}
The resolution of this problem would become feasible if one manages to extend the $L_2$-quotient algorithms to all almost simple groups with socle $PSL(2,q)$. We have no idea however on how easy or difficult the latter is.

\begin{problem}
Enumerate C-group representations of groups with socle $PSL(2,q)$.
\end{problem}

As pointed out earlier, small Ree groups and Suzuki groups are very similar. It would therefore be interesting to try to prove a theorem similar to Theorem~\ref{cgroupssz} for small Ree groups. This is stated in the following problem.
\begin{problem}
Determine the maximal rank of a C-group with group $\Ree(q)$.
\end{problem}

%\subsection{C-group representations of $S_n$}
As in the case of Suzuki groups and groups $PSL(2,q)$ and $PGL(2,q)$, it is possible to obtain classification theorems for the symmetric groups when removing the string condition.
Theorem 1 of~\cite{fl} was obtained as a corollary of the following result
that gives the C-groups of rank $n-1$ for $S_n$ when $n\geq 7$.
It is due to Cameron and Philippe Cara but rephrased here in the framework of C-groups.
\begin{theorem}\label{CaCa}\cite{CaCa:ta}
For $n\geq 7$, $G$ a permutation group of degree $n$ and  $\{\rho_0, \ldots, \rho_{n-2}\}$ a set of involutions of $G$, $\Gamma := (G, \{\rho_0, \ldots, \rho_{n-2}\})$ is a C-group representation of rank $n-1$ of $G$ if and only if the permutation representation graph of $\Gamma$ is a tree $\mathcal{T}$ with $n$ vertices.
Moreover, the $\rho_i$'s are transpositions, $G\cong S_n$ and the Coxeter diagram of $\Gamma$ is the line graph of $\mathcal{T}$.
%Finally, every such C-group gives a regular hypertope of rank $n-1$ for $S_n$.
\end{theorem}

%We refer to~\cite{flw} for the definition of hypertope.

Together with Fernandes, we managed in~\cite{hypsym} to give a classification of C-groups of rank $n-2$ for $S_n$ when $n\geq 9$.
Recall that a 2-transposition is an involution that is the product of two transpositions.

\begin{theorem}\label{hypsym:main}
Let  $n\geq 9$. Let $\{\rho_0, \ldots, \rho_{n-3}\}$ be a set of involutions of $S_n$. Then $\Gamma := (S_n, \{\rho_0, \ldots, \rho_{n-3}\})$ is a C-group representation of rank $n-2$ of $S_n$ if and only if 
its permutation representation graph belongs to one of the following three families, up to a renumbering of the generators, where  $\rho_2, \ldots, \rho_{n-3}$ are transpositions corresponding to the edges of a tree with $n-3$ vertices and the two remaining involutions, $\rho_0$ and $\rho_1$, are either transpositions or $2$-transpositions (with at least one of them being a $2$-transposition).
$$(A)
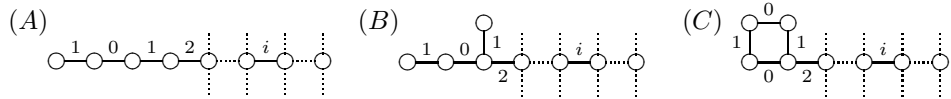
\xymatrix@-1.3pc{&&&&&&&\\*+[o][F]{}   \ar@{-}[r]^1& *+[o][F]{}   \ar@{-}[r]^0 & *+[o][F]{} \ar@{-}[r]^1& *+[o][F]{}\ar@{-}[r] ^2& *+[o][F]{}\ar@{.}[r]\ar@{.}[d]\ar@{.}[u]& *+[o][F]{}\ar@{-}[r]^i\ar@{.}[d]\ar@{.}[u]& *+[o][F]{}\ar@{.}[r] \ar@{.}[d]\ar@{.}[u]&*+[o][F]{}\ar@{.}[d]\ar@{.}[u]\\
&&&&&&&} \quad  (B) \xymatrix@-1.3pc{&&*+[o][F]{}  &&&&&\\*+[o][F]{}   \ar@{-}[r]^1& *+[o][F]{}   \ar@{-}[r]^0 & *+[o][F]{} \ar@{-}[u]_1\ar@{-}[r]_2 & *+[o][F]{}\ar@{.}[r]\ar@{.}[d]\ar@{.}[u]& *+[o][F]{}\ar@{-}[r]^i\ar@{.}[d]\ar@{.}[u]& *+[o][F]{}\ar@{.}[r] \ar@{.}[d]\ar@{.}[u]&*+[o][F]{}\ar@{.}[d]\ar@{.}[u]\\
&&&&&&&} (C) \xymatrix@-1.3pc{ *+[o][F]{}   \ar@{-}[r]^0\ar@{-}[d]_1  & *+[o][F]{} \ar@{-}[d]^1&&&&&& \\
*+[o][F]{}   \ar@{-}[r]_0  & *+[o][F]{} \ar@{-}[r]_2 & *+[o][F]{}\ar@{.}[r]\ar@{.}[d]\ar@{.}[u]& *+[o][F]{}\ar@{-}[r]^i\ar@{.}[d]\ar@{.}[u]& *+[o][F]{}\ar@{.}[r] \ar@{.}[d]\ar@{.}[u]&*+[o][F]{}\ar@{.}[d]\ar@{.}[u]\\
&&&&&&&} $$
%Moreover, every such C-group gives a regular hypertope of rank $n-2$ for $S_n$.
\end{theorem}

\section*{Acknowledgements}
The author would like to thank two anonymous referees whose fruitful comments improved this survey.

\bibliographystyle{amsalpha}

\end{document}